\documentclass[reqno]{amsart}
\usepackage{amssymb, latexsym}
\newtheorem{theorem}{Theorem}[section]

\newtheorem{proposition}[theorem]{Proposition}
\newtheorem{lemma}[theorem]{Lemma}
\newcommand{\p}{\varphi}

\newcommand{\ud}{\mathbb D}

\newcommand{\pdd}{\p:\ud\to\ud}
\newcommand{\JC}{Julia-Carath\'{e}odory\ }
\newcommand{\thmskp}{\vskip12pt}

\begin{document}

\title[Extremal Non-Compactness of Composition Operators]{Extremal Non-compactness of Composition Operators with Linear Fractional Symbol}
\author{Estelle L. Basor*}
\author{Dylan Q. Retsek}
\thanks{*The first author was supported in part  
by National Science Foundation grant DMS-0200167.}
\date{July 19, 2005}

\begin{abstract}
We realize norms of  most composition operators $C_\p$ acting on the Hardy space with linear fractional symbol as roots
of hypergeometric functions. This realization leads to simple necessary and sufficient conditions on $\p$ for $C_\p$ to exhibit
extremal non-compactness, establishes equivalence of cohyponormality and cosubnormality of composition operators with linear
fractional symbol, and yields a complete classification of those linear fractional $\p$ that induce
composition operators whose norms are determined by the action of the adjoint $C_\p^*$ on the normalized reproducing kernels in $H^2$.
\end{abstract}

\maketitle

\section{Introduction}\label{S:intro}

Let $\ud=\{z\in\mathbb{C} : |z|<1\}$ be the unit disk in the complex plane.
The classical Hardy space $H^2$ consists of those holomorphic functions on $\ud$ whose Taylor coefficients about the origin are square summable.
The space $H^2$ is a Hilbert space with inner product given by
\[
\left\langle \sum_{j=0}^\infty a_j z^j , \sum_{k=0}^\infty b_k z^k \right\rangle
= \sum_{j=0}^\infty a_j \bar{b_j}\ .
\]

As a consequence of Cauchy's integral formula, point evaluation is a bounded
linear functional on $H^2$. Thus, for each $w\in\ud$, the Riesz representation theorem
guarantees a unique function $K_w \in H^2$ such that
\[
f(w) = \langle f, K_w \rangle
\]
for all $f\in H^2$. The function $K_w$ is called the reproducing kernel at the point $w$
and is given by $K_w(z) = (1-\bar{w}z)^{-1}$. 

Given analytic $\pdd$, the composition operator $C_\p$ is defined on $H^2$ by $C_\p f = f\circ\p$.
A direct application of Littlewood's subordination principle  \cite{Lit} shows that $C_\p$
is in fact a bounded linear operator on $H^2$. Though boundedness is established with relative
ease, in only special cases is $||C_\p||$ known. In general, calculating the norm of a
composition operator appears to be a difficult problem. 

Recall that a bounded linear operator $T$ on a Hilbert space is \emph{compact}
if the image of the closed unit ball under $T$ has compact closure. Thus, one way to try to compute $||C_\p||$ is to compare it to the essential norm
\[
||C_\p||_e = \inf_{K \text{cpt}}||C_\p+K||.
\]
While it is clear that $||C_\p||_e \leq ||C_\p||$, this inequality may be strict. For example,
if $\p(z) = 2/(3-z)$, then $0 < ||C_\p||_e < ||C_\p||$ (see \cite{BHSF}, Theorem 3.9). Because computing essential norms on the Hardy space is
often easier than computing norms, one would be interested to know precisely which self-maps
$\p$ of the unit disk result in composition operators having norm equal to essential norm. Such operators are said to exhibit
\emph{extremal non-compactness}.

In this paper we obtain necessary and sufficient conditions on the linear-fractional map $\pdd$ for the equality $||C_\p||_e = ||C_\p||$ to hold. Our method involves interpreting the norm representation recently discovered by Bourdon, Fry, Hammond, and Spofford \cite{BHSF} in terms 
of hypergeometric function theory. In the next section we develop the background necessary for the sequel.

\section{Background}\label{S:back}

For any holomorphic map $\pdd$ we have
\begin{equation}\label{E:bigineq}
\frac{1}{1-|\p(0)|^2} \leq ||C_\p||^2 \leq \frac{1+|\p(0)|}{1-|\p(0)|}
\end{equation}
(see, e.g., \cite{CoM}, Corollary 3.7). Equality is achieved on the left when $\p$ is constant
and on the right when $\p$ is inner \cite{Nor}.
The precise value of $||C_\p||$ is known in a few other cases as well. For example,
when $\p(0)=0$, equation (\ref{E:bigineq}) implies that $||C_\p||=1$. Or when $\p(z)=sz+t$ ($s$
and $t$ necessarily satisfying $|s|+|t|\leq 1$), C. Cowen has computed
\[
 ||C_\p||=\sqrt{\frac{2}{1+|s|^2-|t|^2+\sqrt{(1-|s|^2+|t|^2)^2-4|t|^2}}}
\]
in \cite{Cow}.

In attempting to compute the norm of a linear operator $T:\mathcal{H}\to\mathcal{H}$, one often works
with the positive operator $T^*T$ and attempts to compute $||T^*T||=||T||^2$. The ingredients that can
make this approach particularly successful when dealing with the composition operator $C_\p:H^2\to H^2$ are twofold.

First, Hammond shows in \cite{Ham} that if $\p$ is an analytic self-map of $\ud$ such that $||C_\p|| > ||C_\p||_e$, then
there exists a nowhere vanishing function $g\in H^2$ such that $C_\p^*C_\p g = ||C_\p||^2g$. 

Second, when $\p(z)=(az+b)/(cz+d)$ is a non-constant linear fractional map the adjoint of $C_\p$ is computable via Cowen's formula
\[
C_\p^* = T_\gamma C_\sigma T_\nu^*
\]
where $\gamma(z)=1/(-{\bar{b}}z+{\bar{d}})$, $\nu(z)=cz+d$, $\sigma(z)=({\bar{a}}z-{\bar{c}})/(-{\bar{b}}z+{\bar{d}})$ and $T_h$ is the analytic
Toeplitz operator with symbol $h$.

Using these two facts, the quest to compute $||C_\p||$ when $\p(z)=(az+b)/(cz+d)$ and $||C_\p|| > ||C_\p||_e$ is reduced to looking for the largest eigenvalue
of $C_\p^*C_\p$. Toward this end, Hammond established in \cite{Ham} that if $f$ is an eigenvector for $C_\p^*C_\p$ with corresponding eigenvalue
$\lambda$, the functional equation
\begin{equation}\label{E:hamfunc}
\lambda f(z)=\psi(z)f(\tau(z))+\chi(z)\lambda f(0)
\end{equation}
holds at every $z$ for which the auxiliary maps
\begin{equation}\label{E:auxmaps}
\psi(z)=\frac{({\overline{ad}}-{\overline{bc}})z}{({\bar{a}}z-{\bar{c}})(-{\bar{b}}z+{\bar{d}})},\  \chi(z)=\frac{{\bar{c}}}{-{\bar{a}}z+{\bar{c}}},\ \text{and}\ 
\tau(z)=\p(\sigma(z))
\end{equation}
are defined.

Iteration of equation (\ref{E:hamfunc}) ultimately leads to the following theorem of Bourdon, Fry, Hammond and Spofford which appears as theorem 3.5 and Corollary 3.6
in \cite{BHSF}.
\begin{theorem}\label{T:BHSFrep}
Let $\p(z)=(az+b)/(cz+d)$ be a non-automorphic linear fractional mapping that fixes the point 1. Then $||C_\p|| > ||C_\p||_e$ if and only if
there exists a number $\Lambda > ||C_\p||_e^2$ such that
\begin{equation}\label{E:key}
\sum_{k=0}^\infty \chi(\tau^{[k]}(\p(0)))\left[\prod_{m=0}^{k-1}\psi(\tau^{[m]}(\p(0)))\right]\left(\frac{1}{\Lambda}\right)^{k+1}=1
\end{equation}
Moreover, the largest $\Lambda$ (when there is one) for which equation (\ref{E:key}) holds is $||C_\p||^2$.
\end{theorem} 

In general, it is not easy to tell whether equation (\ref{E:key}) has any solutions. In special circumstances, however, one can use this equation to good effect.
For instance, in \cite{BHSF} Bourdon et al. show that if $\p(z)=(\alpha-1)/(\alpha-z)$ for some $\alpha > 1$, then $||C_\p|| > ||C_\p||_e$ (the specific case
$\alpha=3$ was mentioned in the Introduction). The proof is via the intermediate value theorem and depends critically upon the fact that the coefficients of the
series in (\ref{E:key}) are all real for this particular kind of $\p$.

It is a major goal of this paper to circumvent this dependence on real coefficients in the series of Theorem \ref{T:BHSFrep} and show that $||C_\p||_e$ very
rarely equals $||C_\p||$ for linear fractional $\p$. The key observation toward this end is the subject of the next section.

\section{Hypergeometric Series}\label{S:hypergeo}

A hypergeometric series is a series $\sum a_k$ such that the ratio $a_{k+1}/a_k$ is a rational function of $k$. In general,
the numerator and denominator of $a_{k+1}/a_k$ will be polynomials in $k$ of arbitrary degree. Of particular importance to us, however, is the 
special case
\begin{equation}\label{E:F21}
_2F_1(a,b;c;z) \equiv \sum_{k=0}^\infty \frac{(a)_k(b)_k}{(c)_k k!}z^k
\end{equation}
where $(\alpha)_0=1$ and $(\alpha)_k=(\alpha)(\alpha+1)\cdots(\alpha+k-1)$ for $k=1,2,3,\ldots$.

The function $_2F_1(a,b;c;z)$ is called the hypergeometric series of variable $z$ with parameters $a,b,c$. By the ratio test, the series converges
absolutely for $|z|<1$ and therefore defines an analytic function on $\ud$.

Thanks to Euler, Gauss, Riemann and others much is known about the function $_2F_1$, as well as generalized hypergeometric functions having more parameters.
We shall need various facts concerning $_2F_1$ along the way, but for a full treatment see the Bateman manuscript \cite{Bat} or the treatise
of Andrews, Askey and Roy \cite{AAR}.

The remainder of this section is devoted to showing that the series in (\ref{E:key}) is in fact hypergeometric whenever $\p$ is a non-automorphic,
non-affine linear fractional self-map of $\ud$ (extreme non-compactness of $C_\p$ is already settled for automorphic and affine $\p$; see \cite{Sha1},
\cite{Sha2}, \cite{Cow}, and \cite{BHSF}). We begin with several lemmas.

\begin{lemma}\label{L:qdrep}
Suppose $\pdd$ is a non-affine linear fractional map that fixes the point 1. Then $\p$ is of the form
\begin{equation}\label{E:qdrep}
\p(z) = \frac{(1+q+qd)z+(d-q-qd)}{z+d}
\end{equation}
where $q=\p'(1)>0$ and $d\in\mathbb{C}\setminus\overline{\ud}$.
\end{lemma}

\noindent\emph{Proof.} Because $\p$ is not affine, the linear coefficient in the denominator is nonzero. We may therefore
assume without loss of generality that $\p(z)=(az+b)/(z+d)$. Since $\p(1)=1$, $a+b=1+d$. Thus,
\[
\p(z)=\frac{az+(1+d-a)}{z+d}\,.
\]
Differentiating, we find that
\begin{equation}\label{E:phideriv}
\p'(z)=\frac{(a-1)(d+1)}{(z+d)^2}\,.
\end{equation}
Now, because $\p(1)=1$ the \JC theorem implies that $q\equiv \p'(1) > 0$. Substituting $z=1$ in (\ref{E:phideriv}) yields
$a=1+q+qd$ and
\[
\p(z) = \frac{(1+q+qd)z+(d-q-qd)}{z+d}
\]
as claimed.\qed

\thmskp

We note that not all maps of the form in the above lemma are self-maps of the disk. Indeed, $\p(0)=1-q(d+1)/d$ shows that $q$ and $d$ cannot be chosen independently
if one hopes to obtain a self-map. Examination of the conformally equivalent map $f=T\circ\p\circ T^{-1}$, where $T(z)=(1+z)/(1-z)$ maps the disk onto the
right half-plane, shows that $\pdd$ if and only if $\Re\{(d-1)/(d+1)\}\geq q$. This condition will play a role later.

Lemma \ref{L:qdrep} is useful in that it pares down the number of parameters we must deal with. Because we wish to eventually analyze the series
of Theorem \ref{T:BHSFrep}, we would also be well served to simplify the auxiliary maps associated to the linear fractional map $\p$. Recall that if
$\p(z)=(az+b)/(cz+d)$, then $\sigma(z)=({\bar{a}}z-{\bar{c}})/(-{\bar{b}}z+{\bar{d}})$ and $\tau(z)=\p(\sigma(z))$. The next lemma gives an explicit
representation of $\tau$ in terms of the parameters $q$ and $d$ of Lemma \ref{L:qdrep}.

\begin{lemma}\label{L:taurep}
Suppose $\pdd$ is of the form
\[
\p(z) = \frac{(1+q+qd)z+(d-q-qd)}{z+d}\,.
\]
Then the associated map $\tau$ is of the form
\begin{equation}\label{E:taurep}
\tau(z) = \frac{(1-b)z+b}{-bz+1+b}
\end{equation}
where
\begin{equation}\label{E:brep}
b=\frac{|d|^2-q|1+d|^2-1}{q|1+d|^2}\,.
\end{equation}
\end{lemma}

\noindent\emph{Proof.} It is customary (and exceedingly useful) to work with the ``right half-plane version" of $\tau$. Specifically, let
$f=T\circ\tau\circ T^{-1}$ where $T(z)=(1+z)/(1-z)$ is the conformal mapping that takes the unit disk onto the right half-plane. For general
$\p(z)=(az+b)/(cz+d)$, a lengthy and tedious computation shows that
\begin{equation}\label{E:rhptau}
f(z)=T(\tau(T^{-1}(z)))=\frac{Az+B}
                             {Cz+D}
\end{equation}
where
\begin{align}
A &= |a|^2-|b|^2-|c|^2+|d|^2+(c{\bar{a}}-a{\bar{c}})+(b{\bar{d}}-d{\bar{b}})\notag\\
B &=-|a|^2+|b|^2-|c|^2+|d|^2-(c{\bar{a}}+a{\bar{c}})+(b{\bar{d}}+d{\bar{b}})\notag\\
C &=-|a|^2+|b|^2-|c|^2+|d|^2+(c{\bar{a}}+a{\bar{c}})-(b{\bar{d}}+d{\bar{b}})\notag\\
D &= |a|^2-|b|^2-|c|^2+|d|^2-(c{\bar{a}}-a{\bar{c}})-(b{\bar{d}}-d{\bar{b}})\notag
\end{align}

Under the present hypotheses, $a=1+q+qd$, $b=d-q-qd$, and $c=1$. In particular, $C=0$ and computation of $A$, $B$, and $D$ yields
\[
f(z)=\frac{2q|1+d|^2z-4[1+q|1+d|^2-|d|^2]}{2q|1+d|^2}=z+\frac{2(|d|^2-q|1+d|^2-1)}{q|1+d|^2}\equiv z+2b\,.
\]
The utility of this right half-plane version of $\tau$ is now evident; the map $f(z)$ is just translation of the half-plane by the nonnegative real number $2b$.

To complete the proof, we simply unravel $\tau$ and compute
\[
\tau(z)=T^{-1}(f(T (z)))=T^{-1}\left(\frac{1+z}{1-z}+2b\right)=\frac{(1-b)z+b}{-bz+1+b}\,.
\]
\qed

We remark that the form of the coefficients of the series in (\ref{E:key}) indicates that it will be important to know not just $\tau$ itself, but its
$k$-fold iterates $\tau^{[k]}$ as well. As in the proof of Lemma \ref{L:taurep} this is best approached by examining the iterates $f^{[k]}$ of the right
half-plane version of $\tau$. Since it is clear that $f^{[k]}(z)=z+k(2b)$, we immediately obtain
\begin{equation}\label{E:tauiterate}
\tau^{[k]}(z)=T^{-1}(f^{[k]}(T (z)))=\frac{(1-kb)z+kb}{-kbz+1+kb}\,.
\end{equation}

With this firm grip on the auxiliary map $\tau$, we are now in a position to show that the series
\[
\sum_{k=0}^\infty \chi(\tau^{[k]}(\p(0)))\left[\prod_{m=0}^{k-1}\psi(\tau^{[m]}(\p(0)))\right]x^{k+1}
\]
is hypergeometric whenever $\p$ is a non-affine linear fractional self map of the unit disk. 
For convenience we remind the reader of the auxiliary maps corresponding to $\p(z)=(az+b)/(cz+d)$ that were discussed in Section \ref{S:back}:
\[
\psi(z)=\frac{({\overline{ad}}-{\overline{bc}})z}{({\bar{a}}z-{\bar{c}})(-{\bar{b}}z+{\bar{d}})}\ \text{and}\  \chi(z)=\frac{{\bar{c}}}{-{\bar{a}}z+{\bar{c}}}
\]

\begin{theorem}\label{T:hypergeorep}
Let $\pdd$ be a non-affine, non-automorphic linear fractional map that fixes the point 1. Then
\[
\p(z) = \frac{(1+q+qd)z+(d-q-qd)}{z+d}
\]
and 
\begin{equation}\label{E:hypergeorep}
\sum_{k=0}^\infty \chi(\tau^{[k]}(\p(0)))\left[\prod_{m=0}^{k-1}\psi(\tau^{[m]}(\p(0)))\right]x^{k+1}=
1- \,_2F_1(\alpha,\beta;\delta;x/q)
\end{equation}
where
%\begin{align}
%\alpha &= \frac{2+2d}{qb|1+d|^2}\notag\\
%
%
%\beta &=\frac{2d+2d{\bar{d}}}{qb|1+d|^2}\\
%
%
%\delta &= \frac{2}{qb}\notag
%
%
%\end{align}

\[
\alpha = \frac{1+d}{qb|1+d|^2}, \quad \beta =\frac{d+d{\bar{d}}}{qb|1+d|^2}, \quad \text{and}\quad  \delta = \frac{1}{qb}\,.\notag
\]

\end{theorem}

\noindent\emph{Proof.} The form of $\p$ is guaranteed by Lemma \ref{L:qdrep}. Thus, $\p(0)=(d-q-qd)/d$ and the remarks following Lemma \ref{L:taurep}
show that
\[
\tau^{[k]}(\p(0))=\frac{(1-kb)\p(0)+kb}{-kb\p(0)+1+kb}=1-\frac{q(1+d)}{d+qb(1+d)k}
\]
(where $b$ is given by Equation (\ref{E:brep}) above).

Let $\sum_{k=0}^\infty a_kx^{k+1}$ be the series on the left side
of equality (\ref{E:hypergeorep}). To show that the series is hypergeometric, we analyze the
quotient
\begin{equation}\label{E:akfrac}
\frac{a_{k+1}}{a_k} = \frac{\chi(\tau^{[k+1]}(\p(0)))\psi(\tau^{[k]}(\p(0)))}{\chi(\tau^{[k]}(\p(0)))}\,.
\end{equation}

This appears to be a daunting task. Upon substituting the formulas for the functions $\psi$ and $\chi$ (with the coefficients for $\p$ already
in terms of $q$ and $d$) and the value for $\tau^{[k]}(\p(0))$ obtained above, a fair amount of simplification takes place. Continuing from (\ref{E:akfrac})
we have
\begin{align}
\frac{a_{k+1}}{a_k} &= \left\{\frac{q[d(1+{\bar{d}})-q|1+d|^2+qb|1+d|^2k]}{[q(|d|^2-1)+q^2(b-1)|1+d|^2+q^2b|1+d|^2k]}\right\}\notag\\
                    &\times \left\{\frac{d(1+{\bar{d}})+qb|1+d|^2+qb|1+d|^2k]}{[q(|d|^2-1+|1+d|^2)-q^2|1+d|^2+q^2b|1+d|^2k]}\right\}\notag
\end{align}

While we are encouraged because $a_{k+1}/a_k$ indeed appears to be a rational function of $k$, this current incarnation of the quotient is far from practical.
Factoring out a $q$ from each term in the denominator and dividing top and bottom by $(qb|1+d|^2)^2$ we finally arrive at the more palatable
\begin{equation}\label{E:akfracfinal}
\frac{a_{k+1}}{a_k} = \frac{1}{q}\frac{[\alpha+k+1][\beta+k+1]}{[k+2][\delta+k+1]}
\end{equation}
where
\[
\alpha = \frac{1+d}{qb|1+d|^2},\ 
\beta =\frac{d+d{\bar{d}}}{qb|1+d|^2},\ \text{and}\ 
\delta = \frac{1}{qb}
\]
are just as in the statement of the theorem.

It remains to prove that $\sum_{k=0}^\infty a_kx^{k+1}=1-_2F_1(\alpha,\beta;\delta;x/q)$. Consider the series $\sum_{k=0}^\infty b_kx^{k+1}$ where 
the coefficient $b_k$ is given by
\[
b_k=\frac{(\alpha+1)_k(\beta+1)_k}{(\delta+1)_k (k+1)!}\left(\frac{1}{q}\right)^k\,.
\]
By design, the quotient $b_{k+1}/b_k = a_{k+1}/a_k$ satisfies Equation (\ref{E:akfracfinal}). 
Because the first coefficient of any such series determines the remaining coefficients,
we must have
\begin{equation}\label{E:aktobk}
\sum_{k=0}^\infty a_kx^{k+1}=\frac{a_0}{b_0}\sum_{k=0}^\infty b_kx^{k+1}
                            =\frac{a_0}{b_0}\sum_{k=0}^\infty \frac{(\alpha+1)_k(\beta+1)_k}{(\delta+1)_k (k+1)!}\left(\frac{1}{q}\right)^k x^{k+1}\,.
\end{equation}

Now, utilizing the facts that $b_0=1$, $a_0=\chi(\p(0))$, and $\zeta(\zeta+1)_k=(\zeta)_{k+1}$ whenever $\zeta\neq 0$ and $k=0,1,2,\ldots$ (recall
that $\alpha$, $\beta$ and $\delta$ are all nonzero), we continue from (\ref{E:aktobk}) to obtain
\[
q\frac{\chi(\p(0))\delta}{\alpha\beta}
                                \sum_{k=0}^\infty \frac{(\alpha)_{k+1}(\beta)_{k+1}}{(\delta)_{k+1} (k+1)!}\left(\frac{x}{q}\right)^{k+1}
                            =q\frac{\chi(\p(0))\delta}{\alpha\beta}[_2F_1(\alpha,\beta;\delta;x/q)-1]\,.
\]
Finally, a simple computation shows that
\[
\frac{\chi(\p(0))\delta}{\alpha\beta}=-\frac{1}{q}
\]
and the proof is complete.\qed

\thmskp

Because so much is known about hypergeometric functions, we may parlay the information contained in Theorem \ref{T:hypergeorep} to garner information
about the norms of a great many composition operators with linear fractional symbol. Specifically, the work of the next section will show that such
operators are rarely extremally non-compact and, as a consequence, their norms are rarely exhibited by the action of the adjoint $C_\p^*$ acting on the
normalized reproducing kernels in $H^2$.

\section{Extremal Non-compactness}\label{S:noncpt}

Given a linear fractional map $\pdd$, we would like to know the exact conditions under which $||C_\p||=||C_\p||_e$. Because $||\p||_\infty<1$ implies that
$C_\p$ is compact, we limit our discussion to $\p$ for which $||\p||_\infty=1$. In this case, either $\p$ is an automorphism or else $\p(\ud)$ is a disk internally tangent
to $\partial\ud$ at a single point. If $\p$ is an automorphism, then $||C_\p||=||C_\p||_e$ (see, e.g., \cite{Sha2}). We may therefore further limit our discussion
to maps $\p$ such that there exists a single pair of points $\zeta$ and $\eta$ on $\partial\ud$ such that $\p(\zeta)=\eta$.

Having gone this far, we note that without loss of generality we may still further assume that $\zeta=\eta=1$. To wit, Bourdon et al. note in \cite{BHSF} that 
for any unimodular constants $\zeta$ and $\eta$, the operators $C_{\zeta z}$ and $C_{{\bar{\eta}} z}$ are unitary and thus $||C_\p||=||C_\p||_e$ if and only if
$||C_{\zeta z}C_\p C_{{\bar{\eta}} z}||=||C_{\zeta z}C_\p C_{{\bar{\eta}} z}||_e$. But when $\p(\zeta)=\eta$, ${\bar{\eta}}\p(\zeta z)$ is a 
linear fractional self-map of the disk that fixes the point 1. It therefore suffices to know how to compare $||C_\p||$ and $||C_\p||_e$ when $\p(1)=1$.

Finally, affine maps $\p(z)=sz+t$ satisfying $||\p||_\infty=1$ also satisfy $||C_\p||=||C_\p||_e$; Cowen's norm formula in \cite{Cow} and Shapiro's essential norm
formula in \cite{Sha1} may be seen to coincide in this case. Thus, we will assume throughout this section that $\pdd$ \emph{is a non-affine, 
non-automorphic linear fractional map that fixes the point 1}. In the face of the rather restrictive hypotheses of Theorem \ref{T:hypergeorep}, this assumption
lends credibility to the work of Section \ref{S:hypergeo} and, via Lemma \ref{L:qdrep}, allows us to use the representation 
\[
\p(z) = \frac{(1+q+qd)z+(d-q-qd)}{z+d}
\]
freely. We begin with a lemma.

%\newpage

\begin{lemma}\label{L:alphanegint}
Suppose $\pdd$ is given by
\[
\p(z) = \frac{(1+q+qd)z+(d-q-qd)}{z+d}
\]
and the number
\[
\alpha = \frac{1+d}{qb|1+d|^2} \leq -1.
\]
Then $||C_\p||>||C_\p||_e$.
\end{lemma}

\noindent\emph{Proof.} Say $\alpha = -t$ for some $t \geq 1$. Then $d<-1$ (recall that $q$ and $b$ are both positive) and substituting 
the value of $b$ given by (\ref{E:brep}) in terms of $q$ and $d$ we have
\[
-t=\frac{1}{qb(1+d)}=\frac{1+d}{d^2-q(1+d)^2-1}.
\]
Solving for $q$, we obtain $q=(1+t(d-1))/(td+t)$. Because $\p$ is non-automorphic and fixes 1, Shapiro's essential norm formula \cite{Sha1} shows that
$||C_\p||_e^2=1/q$. Thus, since $d+1 \leq (d+1)/t < 0$,
\[
||C_\p||_e^2 = \frac{td+t}{1+t(d-1)}\leq\frac{1+d}{d} = 1+\frac{1}{d} < 1 \leq ||C_\p||^2,
\]
where the last inequality above is true for general self-maps $\p$ by Equation (\ref{E:bigineq}).\qed 

\thmskp

The next proposition shows that the conclusion of the previous lemma holds in general when $d$ is negative.

\begin{proposition}\label{P:negatived}
Suppose $\pdd$ is given by
\[
\p(z) = \frac{(1+q+qd)z+(d-q-qd)}{z+d}.
\]
If $d<-1$, then $||C_\p||>||C_\p||_e$.
\end{proposition}

\noindent\emph{Proof.} By Lemma \ref{L:alphanegint} we may assume that $-1 < \alpha < 0$. According to Theorems \ref{T:BHSFrep} and
\ref{T:hypergeorep},  $||C_\p||>||C_\p||_e$ if and only if the equation $_2F_1(\alpha,\beta;\delta;x/q)=0$ has a solution in the interval $(0,q)$. Thus,
we aim to prove that $_2F_1(\alpha,\beta;\delta;x)=0$ has a root in $(0,1)$ where
\[
\alpha = \frac{1+d}{qb|1+d|^2},\ 
\beta =\frac{d+d{\bar{d}}}{qb|1+d|^2},\ \text{and}\ 
\delta = \frac{1}{qb}\,.
\]

When the number $d$ is real there are some noteworthy relations among the constants $\alpha$, $\beta$, and $\delta$.
First, $\beta=\alpha d$ when $d\in\mathbb{R}$; since $\alpha$ and $d$ have the same sign (again, recall that $b$ and $q$ are positive), we conclude that 
$\beta>0$ whenever $d$ is real. Second, we in general have
\[
\delta-\beta = \frac{1}{qb} - \frac{d+d{\bar{d}}}{qb|1+d|^2} = \frac{1+{\bar{d}}}{qb|1+d|^2}={\bar{\alpha}}\,.
\]
So when $d$ is real, $\delta-\beta=\alpha$.

In general, the series $_2F_1(\alpha,\beta;\delta;x)$ has radius of convergence $R=1$ and diverges at $x=1$ when $\delta=\alpha+\beta$ (terminating
cases $\alpha=-n$ notwithstanding). In this case, the behavior of the series as $x\to 1^{-}$ was determined by Gauss (see \cite{AAR}, Theorem 2.1.3):
\begin{equation}\label{E:gauss}
\lim_{x\to 1^{-}}\frac{_2F_1(\alpha,\beta;\delta;x)}{-\log(1-x)}=\frac{\Gamma(\delta)}{\Gamma(\alpha)\Gamma(\beta)}\,.
\end{equation}
Here $\Gamma(z)=\int_0^\infty t^{z-1}e^{-t}\,dt$ is the Gamma function familiar from complex analysis, analytically continued to the plane
minus zero and the negative integers. 

Since $-1 < \alpha < 0$, $\beta>0$, and $\delta>0$, the right hand side of (\ref{E:gauss}) is a negative
real number. Because $-\log(1-x)\to\infty$ as $x\to1^{-}$, we conclude that $_2F_1(\alpha,\beta;\delta;x)\to -\infty$ as $x\to1^{-}$. The function
$_2F_1(\alpha,\beta;\delta;x)$ is continuous on $0 \leq x < 1$ and $_2F_1(\alpha,\beta;\delta;0)=1$, so the intermediate value theorem implies that
$_2F_1(\alpha,\beta;\delta;x)$ has a root in $(0,1)$.\qed

\thmskp

In spirit, the analysis of the series $_2F_1$ in the proposition above is very much akin to that found in Theorem 3.7 of \cite{BHSF}. In fact,
the astute reader may have already noticed that the map $\p$ of Proposition \ref{P:negatived} is indeed the member $\phi_{r,s}$ of the Cowen-Kriete
family studied there with $r=(d-q-qd)^{-1}$ and $s=q$ (about which we shall have more to say later). In this sense, the hypergeometric approach 
simply offers an alternative light in which to cast these types of series equations. 

The greatest strength of the hypergeometric approach, however, 
is its capacity to give information when the coefficients of the linear fractional map
$\p$ are complex-valued. One readily sees that the proof of the preceeding proposition fails (on several counts) when $d$ is not real. 
The remedy relies on a lovely transformation of Pfaff (see \cite{AAR}, Theorem 2.2.5 for a proof).

\begin{theorem}\label{T:pfaff} (Pfaff's Transformation)
\begin{equation}\label{E:pfaff}
_2F_1(a,b;c;x)=(1-x)^{-a}\, _2F_1\left(a,c-b;c;\frac{x}{x-1}\right)
\end{equation}
\end{theorem}

We remark that one must take care in interpreting the right side of (\ref{E:pfaff}). For $0<x<1$ the argument $x/(x-1)$ takes values in $(-\infty, 0)$, yet
the series definition of $_2F_1$ is valid only in the unit disk. There are, however, various integral representations of hypergeometric functions that
allow analytic continuation to the whole plane minus the ray $[1,\infty)$. Interpreted this way, Pfaff's transformation is valid and we henceforth let 
context dictate what we mean by the function $_2F_1$. For a much more extensive treatment of analytic continuation of 
hypergeometric functions refer to \cite{AAR} or \cite{Bat}. 

We are now in a position to prove the main theorem of this section.

\begin{theorem}\label{T:main}
Suppose $\pdd$ is given by
\[
\p(z) = \frac{(1+q+qd)z+(d-q-qd)}{z+d}.
\]
Then $||C_\p||=||C_\p||_e$ if and only if $d>1$.
\end{theorem}

\noindent\emph{Proof.} We know already by Theorems \ref{T:BHSFrep} and \ref{T:hypergeorep} that $||C_\p||>||C_\p||_e$ if and only if the series
 $_2F_1(\alpha,\beta;\delta;x)$ has a root in $(0,1)$.

Suppose $d>0$. The coefficients of the series $_2F_1(\alpha,\beta;\delta;x)$ are then clearly positive. Since
$_2F_1(\alpha,\beta;\delta;0)=1$, the series has no roots in $(0,1)$ and therefore $||C_\p||=||C_\p||_e$.

For the converse, suppose that $d$ is not positive. Proposition \ref{P:negatived} handles the case $d<-1$, so we are left to consider whether
$_2F_1(\alpha,\beta;\delta;x)=0$ holds for any $x$ in $(0,1)$ when $0 < |\arg d| < \pi$. Because $\delta-\beta={\bar{\alpha}}$, Pfaff's transformation implies
\begin{equation}\label{E:pfaffrep}
_2F_1(\alpha,\beta;\delta;x)=(1-x)^{-\alpha}\, _2F_1\left(\alpha,{\bar{\alpha}};\delta;\frac{x}{x-1}\right)\,.
\end{equation}

The utility of equation (\ref{E:pfaffrep}) is that the hypergeometric series on the right side is \emph{real} for all $x\in(0,1)$. This is easy to see when
$0 < x < 1/2$; in this case, $x/(x-1)$ is in the unit disk and the coefficients in the series definition of $_2F_1$ are all real. When $1/2 \leq x < 1$, things
are less clear. For these values of $x$, the argument $x/(x-1)$ lies in the interval $(-\infty, -1]$ and we must take $_2F_1$ to mean the analytic continuation
of the series valid in the disk. Letting $t=x/(x-1)$ in this case, the identity
\begin{align}
\frac{_2F_1(\alpha,{\bar{\alpha}};\delta;t)}{\Gamma(\delta)}
           &= \frac{\Gamma({\bar{\alpha}}-\alpha)}{\Gamma({\bar{\alpha}})\Gamma(\delta-\alpha)}(-t)^{-\alpha}\,
                       _2F_1(\alpha, 1-\delta+\alpha; 1-{\bar{\alpha}}+\alpha;1/t)\\
           &+ \frac{\Gamma(\alpha-{\bar{\alpha}})}{\Gamma(\alpha)\Gamma(\delta-{\bar{\alpha}})}(-t)^{-{\bar{\alpha}}}\,
                       _2F_1({\bar{\alpha}}, 1-\delta+{\bar{\alpha}}; 1-\alpha+{\bar{\alpha}};1/t),\notag
\end{align}
which follows from Barnes' integral representation for $_2F_1$ and is valid provided $\alpha-{\bar{\alpha}}$ is not an integer (see \cite{Bat}, Section
2.1.4), shows that 
\[
{\overline{_2F_1(\alpha,{\bar{\alpha}};\delta;x/(x-1))}} =\, _2F_1(\alpha,{\bar{\alpha}};\delta;x/(x-1))
\]
holds for $1/2 \leq x < 1$ (due to the symmetry in $\alpha$ and ${\bar{\alpha}}$).

The point here is that $(1-x)^{-\alpha}$ is nonzero for $0 < x < 1$, so according to Equation (\ref{E:pfaffrep}) $_2F_1(\alpha,\beta;\delta;x)$ has 
a root in $(0,1)$ if and only if  $_2F_1(\alpha,{\bar{\alpha}};\delta;t)$ has a root $t$ in $(-\infty,0)$. Since we have established that
$_2F_1(\alpha,{\bar{\alpha}};\delta;t)$ is real for these values of $t$, we may once again employ the intermediate value theorem to guarantee the desired root.
Specifically, it now suffices to show that $_2F_1(\alpha,{\bar{\alpha}};\delta;t)<0$ for some $t\in(-\infty, 0)$.

Toward this goal, observe that Equation (28) may be written
\begin{align}
_2F_1(\alpha,{\bar{\alpha}};\delta;t)
           &= \frac{\Gamma({\bar{\alpha}}-\alpha)\Gamma(\delta)}{\Gamma({\bar{\alpha}})\Gamma(\delta-\alpha)(-t)^{\alpha}}\,
                       \left[1+O\left(\frac{1}{t}\right)\right]\notag\\
           &+ \frac{\Gamma(\alpha-{\bar{\alpha}})\Gamma(\delta)}{\Gamma(\alpha)\Gamma(\delta-{\bar{\alpha}})(-t)^{{\bar{\alpha}}}}\,
                       \left[1+O\left(\frac{1}{t}\right)\right]\quad (t\to-\infty) \notag
\end{align}
Letting $s=-t$, it therefore suffices to show that
\[
(A+Bi)s^{-\alpha}+(A-Bi)s^{-{\bar{\alpha}}} < 0
\]
for arbitrarily large positive values of $s$ and
\[
A+Bi \equiv  \frac{\Gamma({\bar{\alpha}}-\alpha)\Gamma(\delta)}{\Gamma({\bar{\alpha}})\Gamma(\delta-\alpha)}\,.
\]
Now, if $\alpha=x+iy$, then 
\[
(A+Bi)s^{-\alpha}+(A-Bi)s^{-{\bar{\alpha}}} = 2s^{-x}[A\cos(y\ln(s))+B\sin(y\ln(s))]\,.
\]
Note that $A$ and $B$ cannot both be zero since the function $_2F_1(\alpha,{\bar{\alpha}};\delta;t)$ is not identically zero for $t<0$.
Also note that because $0 < |\arg d| < \pi$, $\alpha$ is not real. In particular, $y\neq 0$ and, consequently, the function
\[
G(\omega)=A\cos(y\omega)+B\sin(y\omega)
\]
is negative for arbitrarily large positive values of $\omega$. For instance, if $A>0$ then
$G((2k+1)\pi/y)=-A$ for each $k\in\mathbb{N}$; the possibilities $A<0$, $B>0$, and $B<0$ are handled in a similar fashion.\qed

\thmskp

This completes the classification of the extremely non-compact composition operators induced by linear fractional self maps of the disk. To summarize,
there are three possibilities for a linear fractional map $\pdd$:
\begin{itemize}
\item $||\p||_\infty<1$, in which case $||C_\p||>||C_\p||_e=0$.
\item $\p$ is an automorphism, in which case $$||C_\p||=||C_\p||_e=[(1+|\p(0)|)/(1-|\p(0)|)]^{1/2}.$$
\item The image $\p(\ud)$ is a disk internally tangent to $\partial\ud$ at a single point. In this case, either $\p(z)=sz+t$ and 
$||C_\p||=||C_\p||_e$ or else $\p$ is a non-affine, non-automorphic map. If the latter is the case, then $\p$ may be pre- and post-composed
with rotations to produce a non-affine, non-automorphic map $\tilde\p$ that fixes the point 1 and induces a composition operator with the
same norm and essential norm as $C_\p$. This new map must be of the form
\[
\tilde\p(z) = \frac{(1+q+qd)z+(d-q-qd)}{z+d}
\]
and has $||C_{\tilde\p}||=||C_{\tilde\p}||_e$ if and only if $d>1$.
\end{itemize}

Christopher Hammond has observed that this classification settles another closely related operator-theoretic
question. This is the subject of the next section.

\section{Cohyponormality}\label{S:cohypo}

Recall that for a self-map $\p$, the Denjoy-Wolff point of $\p$ is the unique point
$a\in\overline\ud$ such that $|\p'(a)|\leq 1$. If $|a|<1$, then the spectral radius
$r(C_\p)$ of $C_\p$ is 1; if $|a|=1$, then $r(C_\p) = [\p'(a)]^{-1/2}$ (see \cite{CoM}, Th. 3.9).
When working with norms, it is thus sometimes helpful to deal with composition operators that enjoy equality of norm and spectral radius. 
The composition operators exhibiting some form of normality form just such a class (recall that $T:\mathcal{H}\to\mathcal{H}$ is normal if $TT^*=T^*T$).

A Hilbert space operator $T:\mathcal{H}\to\mathcal{H}$ is called \emph{hyponormal} if
$T^*T-TT^* \geq 0$ (equivalently, $||Th|| \geq ||T^*h||$ for all $h\in\mathcal{H}$). The operator $T$ is \emph{subnormal} if 
there exists an auxiliary Hilbert space $\mathcal{K}\supseteq\mathcal{H}$ and a normal operator $N:\mathcal{K}\to\mathcal{K}$ such that $T$ is 
the restriction of $N$ to $\mathcal{H}$. In this section we shall be concerned with deciding when the adjoint $T^*$ satisfies these definitions, in which case the
operator $T$ is said to be \emph{co}hyponormal or \emph{co}subnormal.

The operator $C_\p$ itself is known to be subnormal in just a few cases.
These are: when $\p(z)=cz$ (these are precisely the normal $C_\p$), when $\p$ is
inner and fixes the origin (these are precisely the isometric $C_\p$), and
when $\p$ is linear fractional and $C_\p$ is unitarily equivalent to
$1\oplus C_\psi^*$ with $C_\psi^*$ subnormal.

Perhaps the best theorem on subnormality of $C_\p^*$ is the following result of C. Cowen and
T. Kriete in their article \cite{CoK}.
\begin{theorem}\label{T:cksub}
For $0<s<1$ and $0\leq r \leq 1$, if
\[
\p_{r,s}(z)=\frac{(r+s)z + (1-s)}{r(1-s)z + (1+sr)}\ ,
\]
then $C_{\p_{r,s}}^*$ is subnormal.
\end{theorem}

\thmskp

A nice discussion of the proof may be found in \cite{CoM},
p. 314-17. Under some additional smoothness hypotheses on $\p$, the converse
to Theorem \ref{T:cksub} is true as well \cite{CoK}. The authors there
suggest that the smoothness assumption serves only to facilitate the proof
of the converse and that probably it is true without any additional 
requirements on $\p$.

Using the definitions above, it is not hard to show that any subnormal operator is hyponormal. Though not true in general, for adjoints of composition operators with linear fractional symbol the converse holds as well.

\begin{theorem}\label{T:cohypo}
Let $\pdd$ be linear fractional. If $C_\p$ is cohyponormal, then $C_\p$ is cosubnormal.
\end{theorem}

\noindent\emph{Proof.} Suppose $C_\p^*$ is hyponormal. If $C_\p^*$ is normal, then it's certainly subnormal. We may therefore
assume without loss of generality that $C_\p^*$, hence $C_\p$, is not normal.

Theorem 8.4 of \cite{CoM} then implies that $\p$ has Denjoy-Wolff point $a\in\partial\ud$ with $\p'(a)<1$. Since both
hypo- and subnormality are unitary invariants, we may assume (as usual) that $a=1$. The remainder of the proof proceeds in cases.

If $\p$ is affine, then $\p(z)=sz+(1-s)$ where $0 < s=\p'(1) < 1$. In this case, $\p=\p_{0,s}$ and $C_\p^*$ is subnormal by
Theorem \ref{T:cksub}.

If $\p$ is an automorphism, then a priori it is either elliptic, parabolic, or hyperbolic. Because $\p(1)=1$, $\p$ is not elliptic.
In addition, $\p$ is not parabolic because in this case we have $r(C_\p)=1<||C_\p||$ by Theorem 7.5 of \cite{CoM} (recall that our
cohyponormality assumption guarantees norm and spectral radius coincide). It is possible, however, for $\p$ to be hyperbolic.

Suppose that $\p(z)=\lambda(w-z)/(1-{\bar{w}}z)$ is hyperbolic. In this case, Theorem 7.4 of \cite{CoM} asserts that $r(C_\p)=(\p'(1))^{-1/2}$. Equating
this spectral radius to $||C_\p||$ forces $-1 < w < 0$ and one finds that in actuality $\p=\p_{r,s}$ where $r=1$ and $s=||C_\p||^{-2}$. Thus, when $\p$ is
hyperbolic, Theorem \ref{T:cksub} again implies that $C_\p^*$ is subnormal.

Finally, if $\pdd$ is a non-affine, non-automorphic linear fractional map that fixes the point 1, we are squarely in the realm of 
Theorem \ref{T:main}. That result, combined with the fact that essential norm equals spectral radius for a composition operator with non-automorphic linear fractional
symbol having Denjoy-Wolff point on $\partial\ud$ (see Section 2 of \cite{BHSF} for the details), ensures that $\p$ is of the form
\[
\p(z) = \frac{(1+q+qd)z+(d-q-qd)}{z+d}
\]
for some $d>1$ and $0 < q = \p'(1) < 1$. Set $r=(d-q-qd)^{-1}$ and $s=q$. Since $\Re\{(d-1)/(d+1)\} \geq q$ is necessary for $\p$ to map the disk into
itself, we obtain $0 < r \leq 1$ and $0 < s < 1$. Thus, $\p=\p_{r,s}$ is yet again a member of the Cowen-Kriete family of Theorem \ref{T:cksub} and consequently
$C_\p^*$ is subnormal in this final case as well. $\qed$

\thmskp

We remark that the family of self-maps $\p$ with ``$qd$-representation"
\[
\p(z) = \frac{(1+q+qd)z+(d-q-qd)}{z+d}
\]
where $q>0$ and $|d|>1$ is distinct from the ``extended" Cowen-Kriete family
\[
\p_{r,s}(z)=\frac{(r+s)z + (1-s)}{r(1-s)z + (1+sr)}\ ,
\]
where $0 \leq |r| \leq 1$ and $0 < s < 1$. 

Lemma \ref{L:qdrep} shows that the non-affine members of the extended Cowen-Kriete family certainly have $qd$-representation, but
the self-map $\p(z)=1/(2-z)$, which has $q=1$ and $d=-2$, is not equal to $\p_{r,s}$ for any $r$ and $s$. What has been most important
to us, however, is that self-maps with $qd$-representation, $d>1$, \emph{do} belong the Cowen-Kriete family and consequently enjoy the
attendant nice properties of maps therein.

We should also point out that Theorem \ref{T:cohypo} marks the culmination of steady progress in that direction. Subsuming the recent completion of the story for real $r$ contained in the work of
Bourdon et al. in \cite{BHSF}, we now know that symbols belonging to the extended Cowen-Kriete family induce cohyponormal composition operators if and only if
$0 \leq r \leq 1$.

In the final section of the paper, we consider a detailed example and establish connections between extremal non-compactness and another heavily 
studied operator-theoretic quantity that is often compared to $||C_\p||$.

\section{The Quantity $S_\p^*$ and Geometric Function Theory}\label{S:geotheory}

Another way to try to compute $||C_\p||$ that has drawn a fair amount of attention recently is to compare it to the quantity
\[
S_\p^* := \sup_{w\in\ud} \left\| C_\p^*\left(\frac{K_w}{||K_w||}\right)\right\|
\]
where $K_w$ represents the reproducing kernel at $w$ in $H^2$. While it is clear that $||C_\p|| \geq S_\p^*$, this inequality may be strict. For example,
if $\p(z) = 2/(3-z)$, then $||C_\p|| > S_\p^*$ (\cite{ABT}, Theorem 4.2). On the other hand,
if $\p(z) = sz+t$ for constants satisfying $|s|+|t|\leq 1$, then $||C_\p||=S_\p^*$
(see \cite{Cow}, Proof of Theorem 3, for the case $|s|+|t|=1$; see \cite{BoR}, Proof of Theorem
3.2, for the case  $|s|+|t|<1$).

In \cite{BoR}, the second author and Paul Bourdon show that for $\p$ with $\p(0)\neq 0$ inducing
compact $C_\p$, $||C_\p||=S_\p^*$ if and only if $\p(z) = sz+t$.
As a corollary, non-affine maps $\p$ inducing $C_\p$ with $||C_\p|| = S_\p^*$  are shown to satisfy 
$||C_\p||=||C_\p||_e$ .

As evidenced by Joel Shapiro's celebrated formula for the essential norm of a composition operator on $H^2$ (see \cite{Sha1})
\[
||C_\p||_e^2=\limsup_{|w|\to 1^{-}}\frac{N_\p(w)}{-\log |w|}\ ,
\] 
the compactness of $C_\p$ is intimately related to the way in which $\p$ approaches the 
boundary of the disk, both how fast and how often. Using this formula, it is not difficult
to show that for univalent $\p$, $S_\p^*\geq ||C_\p||_e$. Now, since affine maps induce 
$C_\p$ satisfying either $||C_\p||_e = 0$ or $||C_\p||_e=||C_\p||$, the preceeding results combine to yield the following fact.

\thmskp

\noindent{\emph{Fact}:} If $\p$ is univalent, $\p(0)\neq 0$, and $C_\p$ is
non-compact, then $||C_\p|| = S_\p^*$ if and only if
$||C_\p||=||C_\p||_e$.

\thmskp

The complete classification of extremal non-compactness of composition operators induced by linear fractional maps provided in the previous section
therefore completes the classification of such composition operators satisfying $||C_\p|| = S_\p^*$ as well. Theorem \ref{T:main} adds to the evidence
provided in \cite{BoR} that the set of normalized reproducing kernels is rarely robust enough to exhibit the norm of a composition operator under the 
action of the adjoint $C_\p^*$.

In the case of linear fractional $\p$, there is a particularly nice geometric interpretation of the rarity of this phenomenon.
Let us refer to a linear fractional map $\p$, with image disk $\Delta=\p(\ud)$ tangent at $\eta\in\partial\ud$, as \emph{fast} if for some (hence, all)
$r\in (0, 1)$ we have
\[
d(\partial\Delta, \p(\Gamma_r)\,)=d(\eta, \p(\Gamma_r)\,)\ ,
\]
where $\Gamma_r = \{z: |z|=r\}$ and $d$ is Euclidean distance.

Because the sets $\Gamma_r$ are left invariant under rotations, the same reasoning found at the beginning of Section \ref{S:noncpt} allows us to focus
our attention on maps $\p$ given by
\[
\p(z) = \frac{(1+q+qd)z+(d-q-qd)}{z+d}
\]
In this case, the condition $d>0$ is precisely the condition that $\p$ be fast. Thus, for non-affine, non-automorphic $\p$ that fix the point 1,
$||C_\p|| = S_\p^*$ if and only if $||C_\p||=||C_\p||_e$ if and only if $\p$ is fast. This is intuitively appealing because fast linear fractional
maps approach the boundary of the disk ``as fast as possible" and so should be ``as non-compact as possible".

We conclude with an example. Let $\p(z)=2/(3-z)$. This non-affine, non-automorphic linear fractional self-map of the disk already fixes the point 1 and so has
``$qd-$representation" 
\[
\p(z) = \frac{0z+(-2)}{z+(-3)}
\]
As such, $d=-3$ and we know that $||C_\p|| > ||C_\p||_e$ and $||C_\p|| > S_\p^*$. This is not new information; dating back to \cite{ABT}, this particular map
$\p$ has been extensively studied.

What is more interesting is to examine the family $\{\p_\theta\}_{0\leq\theta\leq2\pi}$ of linear fractional maps given by
\[
\p_\theta(z) = -[\p(\lambda z)-3/4]e^{i\theta}+3/4\,.
\]
Here the unimodular constant $\lambda=\p^{-1}((3-e^{-i\theta})/4)$ is a rotational parameter inserted to guarantee that $\p_\theta$ fixes the point 1 for every 
$\theta$.

We note that $\p_0$ is just the original map $\p$ and that each $\p_\theta$ maps onto the same image disk of radius 1/4 centered at 3/4. But the \emph{way} in which
each $\p_\theta$ maps onto this disk depends on $\theta$. For instance, when $\theta=\pi$ the map
\[
\p_\pi(z) = \frac{(3/2)z+5/2}{z+3}
\]
is fast. Thus, since $d=3$ and $q=1/8$, we have $||C_\p||^2=||C_\p||_e^2=(S_{\p_\pi}^*)^2=8$.

On the other hand, the map
\[
\p_{\frac{\pi}{2}}(z) = \frac{(15+15i)z+(-31+33i)}{20z+(-36+48i)}
\]
is not fast but is ``faster" than the original map $\p(z)=2/(3-z)$. In this case, we use the hypergeometric representation
$_2F_1((-2+6i)/15, (6+2i)/5; 16/15; x)$ to numerically compute
\[
||C_{\p_{\frac{\pi}{2}}}||^2 \approx 3.3764 > 16/5 = ||C_{\p_{\frac{\pi}{2}}}||_e^2
\]

In general,
\[
||C_{\p_\theta}||_e^2=\frac{16}{5+3\cos\theta}
\]
and the closer $\theta$ is to $\pi$, the closer the quotient $Q(\theta)=||C_{\p_\theta}||/||C_{\p_\theta}||_e$ is to 1. We leave the reader with this final
challenge, which to our knowledge has not been met.

\thmskp

\noindent{\emph{Challenge}:} Write the function $Q$ in closed form as a (probably) transcendental function of $\theta$.

\vspace{.2in}
\noindent Department of Mathematics\\
California Polytechnic State University\\
San Luis Obispo, CA 93407, U.S.A.\\
ebasor@calpoly.edu

\vspace{.2in}
\noindent Department of Mathematics\\
California Polytechnic State University\\
San Luis Obispo, CA 93407, U.S.A.\\
dretsek@calpoly.edu

\end{document}